\documentclass[10pt,reqno]{amsart}
\usepackage{amsmath,amsfonts, amssymb, amsthm}
  \usepackage{paralist}
  \usepackage{graphics} 
  \usepackage{epsfig} 
\usepackage{graphicx}  \usepackage{epstopdf}
 \usepackage[colorlinks=true]{hyperref}
\hypersetup{urlcolor=blue, citecolor=red}

\newtheorem{theorem}{Theorem}[section]
\newtheorem{corollary}{Corollary}[section]

\newtheorem{lemma}[theorem]{Lemma}

\newtheorem{conjecture}{Conjecture}[section]

\theoremstyle{definition}

\newtheorem{remark}{Remark}[section]

\def\pmod #1{\ ({\rm{mod}}\ #1)}
\def\Z{\Bbb Z}
\def\N{\Bbb N}

\def\Q{\Bbb Q}

\def\l{\left}
\def\r{\right}
\def\bg{\bigg}
\def\({\bg(}
\def\){\bg)}
\def\t{\text}
\def\f{\frac}

\def\ls{\leqslant}
\def\gs{\geq}

\def\bi{\binom}

\def\eq{\equiv}

\def\da{\delta}

\def\Proof{\noindent{\it Proof}}

\begin{document}
\hbox{Frontiers in Combinatorics and Number Theory 2 (2026), 1--22.}
\medskip

\title[A family of polynomials and related congruences and series]
      {A family of polynomials \\ and related congruences and series}
\author[Zhi-Wei Sun]{Zhi-Wei Sun}


\address{School of Mathematics, Nanjing
University, Nanjing 210093, People's Republic of China}
\email{zwsun@nju.edu.cn}

\keywords{Binomial coefficients, congruence, polynomial, $q$-log-convexity, series for $1/\pi$.
\newline \indent 2020 {\it Mathematics Subject Classification}. Primary 11B65, 11A07; Secondary 05A19.
\newline \indent Supported by the Natural Science Foundation of China (grant no. 12371004).}

\begin{abstract}
In this paper we study a family of polynomials
$$S_n^{(m)}(x):=\sum_{i,j=0}^n\binom ni^m\binom nj^m\binom{i+j}ix^{i+j}\ \ (m,n=0,1,2,\ldots).$$
For example, we show that
$$\sum_{k=0}^{p-1}S_k^{(0)}(x)\equiv\frac x{2x-1}\left(1+\left(\frac{1-4x^2}p\right)\right)\pmod p $$
for any odd prime $p$ and integer $x\not\equiv1/2\pmod p$, where $(\f{\cdot}p)$ denotes the Legendre symbol.

We also formulate some open conjectures on related congruences and series for $1/\pi$. For example, we conjecture that
$$\sum_{k=0}^\infty(7k+1)\frac{S_k^{(2)}(1/11)}{9^k}=\frac{5445}{104\sqrt{39}\,\pi}$$
and
$$\sum_{k=0}^\infty(1365k+181)\frac{S_k^{(2)}(1/18)}{16^k}=\f{1377}{\sqrt2\,\pi}.$$
\end{abstract}
\maketitle

\section{Introduction}

The Ap\'ery numbers given by
$$A_n=\sum_{k=0}^n\bi nk^2\bi{n+k}k^2\ \ (n\in\N=\{0,1,2,\ldots\})$$
play important roles in R. Ap\'ery's proof of the irrationality of $\zeta(3)$ (cf. \cite{A,P}).
The author \cite{S12} conjectured that for any odd prime $p$ we have
$$\sum_{k=0}^{p-1}A_k\eq\begin{cases}4x^2-2p\pmod{p^2}&\t{if}\ p=x^2+2y^2\ (x,y\in\Z),
\\0\pmod{p^2}&\t{if}\ p\eq 5,7\pmod{8}.\end{cases}$$
This was confirmed by C. Wang and the author \cite{WS} via the $p$-adic Gamma function.

The author \cite{S12} introduced the Ap\'ery polynomials
$$A_n(x)=\sum_{k=0}^n\bi nk^2\bi{n+k}k^2x^k\ \ (n\in\N),$$
and investigated their arithmetic properties. For example, Sun \cite{S12} proved that
$$\f1n\sum_{k=0}^{n-1}(2k+1)A_k(x)\in\Z[x]\quad \t{for all}\ n\in\Z^+=\{1,2,3,\ldots\}.$$
V.J.W. Guo and J. Zeng \cite{GZ} confirmed a conjecture of Sun \cite{S12} which states that
$$\f1n\sum_{k=0}^{n-1}(2k+1)(-1)^kA_k(x)\in\Z[x]\quad \t{for all}\ n\in\Z^+.$$

For $m,n\in\N$, we introduce the polynomial
\begin{equation}S_n^{(m)}(x)=\sum_{i=0}^n\sum_{j=0}^n\bi ni^m\bi nj^m\bi{i+j}ix^{i+j}.
\end{equation}
This is motivated by A. Labelle's conjecture (cf. \cite{Q489184})
that
\begin{equation} S_n^{(2)}(1)=A_n\ \ \t{for all}\ n\in\N,
\end{equation}
which has been confirmed by H. Rosengren, M. Alekseyev and A. Labelle in three different ways
(see the answers in \cite{Q489184}).

Our first result is as follows.

\begin{theorem}\label{Th1.1} Let $n\in\N$.

{\rm (i)} For any $m\in\Z^+$, we have the identity
\begin{equation}\label{general}
 S_n^{(m)}(x) =\sum_{k=0}^n\bi nk^2 x^{2k}\(\sum_{j=0}^{n-k}\bi n{j+k}^{m-1}\bi{n-k}jx^j\)^2.
\end{equation}
Also,
\begin{equation}\label{S0}S_n^{(0)}(x)=\sum_{k=0}^{n}x^{2k}\(\sum_{j=0}^{n-k}\bi{j+k}kx^j\)^2.
\end{equation}

{\rm (ii)} We have the identities
\begin{equation}\label{T}S_n^{(1)}(x)=\sum_{k=0}^n\bi nk^2 x^{2k}(1+x)^{2(n-k)}
=\sum_{k=0}^n\bi nk\bi{n+k}kx^{2k}(2x+1)^{n-k}
\end{equation}
and the recurrence
\begin{equation}\label{rec}(n+2)S_{n+2}^{(1)}(x)=(2n+3)(2x^2+2x+1)S_{n+1}^{(1)}(x)-(n+1)(2x+1)^2S_n^{(1)}(x).
\end{equation}

\end{theorem}
\begin{remark} By the Chu-Vandemonde identity (cf. \cite[p.\,22, (3.1)]{G})
and Theorem \ref{Th1.1}(i),  it is easy to see that
$$S_n^{(2)}(1) =A_n\ \ \t{and}\ \ S_n^{(0)}(1)=\bi{2(n+1)}{n+1}-1.$$
\end{remark}

A sequence $(P_n(q))_{n\gs0}$ of polynomials with integer coefficients is said to be {\it $q$-log-convex} if for each positive integer $n$ all the coefficients of the polynomial $$P_{n-1}(q)P_{n+1}(q)-P_n(q)^2\in\mathbb Z[q]$$ are nonnegative.
In 2010, W.Y.C. Chen, R.L. Tang, L.X.W. Wang and A.L.B. Yang \cite{CTWY} proved that
the sequence $(\sum_{k=0}^n\bi nk^2 q^k)_{n\gs0}$ is $q$-log-convex, which was a previous conjecture
of L.L. Liu and Y. Wang \cite{LW}. This, together with \eqref{T}, implies that the sequence $(S_n^{(1)}(q))_{n\gs0}$ is $q$-log-convex. For $\beta_n(q)=\sum_{k=0}^n\bi nk^2\bi{n+k}kq^k$,
the author's conjecture (cf. \cite[Conjecture 4.7]{S14bnk}) on the $q$-log-convexity
 of the sequence $(\beta_n(q))_{n\ge0}$
remains open.

Based on our computation, we formulate the following conjecture.

\begin{conjecture} {\rm (i)} The sequence $(\sum_{k=0}^n\bi nk^3 q^k)_{n\gs8}$ is $q$-log-convex.
Also, for each integer $m\gs2$, the sequence $(\sum_{k=0}^n\bi{n+k}k^mq^k)_{n\gs0}$ is $q$-log-convex.

{\rm (ii)} The sequences $(S_n^{(2)}(q))_{n\gs0}$ and $(S_n^{(3)}(q))_{n\gs2}$
are both $q$-log-convex.
\end{conjecture}

Now we state our second result.

\begin{theorem} \label{Th1.2} Let $p$ be any odd prime, and let $x\in\Z$.
If $x\not\eq1/2\pmod p$, then
\begin{equation}\label{2x-1}\sum_{k=0}^{p-1}S_k^{(0)}(x)
\eq\f x{2x-1}\l(1+\l(\f{1-4x^2}p\r)\r)\pmod p,
\end{equation}
where $(\f{\cdot}p)$ denotes the Legendre symbol.
When $x\eq1/2\pmod p$, we have
\begin{equation}\sum_{k=0}^{p-1}S_k^{(0)}(x)\eq-\da_{p,3}\pmod p.
\end{equation}
\end{theorem}

Theorem \ref{Th1.2} with $x\in\{(p-1)/2,2\}$ yields the following corollary.
\begin{corollary} Let $p$ be an odd prime. Then
\begin{equation}\sum_{k=0}^{p-1}S_k^{(0)}\l(-\f12\r)\eq\f14\pmod p.
\end{equation}
If $p>3$, then we have
\begin{equation}\sum_{k=0}^{p-1}S_k^{(0)}(2)
\eq\f 2{3}\l(1+\l(\f{p}3\r)\l(\f{p}5\r)\r)\pmod p.
\end{equation}
\end{corollary}

Now we give our last theorem.

\begin{theorem} \label{Th1.3} {\rm (i)} For any $n\in\Z^+$, we have
\begin{equation}\label{/n}\f1n\sum_{k=0}^{n-1}S_k^{(1)}(x)\in\Z[x(x+1)]
\end{equation}
and
\begin{equation}\label{6/n}\f{(6,n)}n\sum_{k=0}^{n-1}kS_k^{(1)}(x)\in\Z[x(x+1)],
\end{equation}
where $(6,n)$ is the greatest common divisor of $6$ and $n$.

{\rm (ii)} Let $p$ be any odd prime. Then
\begin{equation}\label{p^2}\f1p\sum_{k=0}^{p-1}S_k^{(1)}(x)\eq 1-(x^{p-1}-1)((x+1)^{p-1}-1)\pmod {p\Z_p[x]},
\end{equation}
where $\Z_p$ is the ring of $p$-adic integers.
\end{theorem}

The classical Ramanujan-type series for $1/\pi$ (cf. \cite{R}) have the form
$$\sum_{k=0}^\infty(ak+b)\f{c_k}{m^k}=\f{\sqrt d}{\pi},$$
where $a,b$ and $m\not=0$ are integers, $d$ is a positive rational number, and $c_k$ (with $k\in\N$) is one of the following products:
$$\bi{2k}k^3,\ \bi{2k}k^2\bi{3k}k,\ \bi{2k}k^2\bi{4k}{2k},\ \bi{2k}k\bi{3k}k\bi{6k}{3k}.$$
One may consult S. Cooper \cite[Chapter 14]{Co} for an introduction to Ramanujan-type series.
In 1997 van Hamme \cite{vH} realized that classical Ramanujan-type serie have corresponding $p$-adic congruences.

The Ap\'ery numbers are related to series for $1/\pi$. In 2002 T. Sato announced the identity
$$\sum_{k=0}^\infty(20k+10-3\sqrt5)\f{A_k}{((\sqrt5+1)/2)^{12k}}=\f{20\sqrt3+9\sqrt{15}}{6\pi}.$$
Motivated by this and the fact that $S_n^{(2)}(1)=A_n$, we seek new series for $1/\pi$ in the form
$$\sum_{k=0}^\infty (ak+b)\f{S_k^{(2)}(c)}{m^k}=\f{\sqrt d}{\pi},$$
where $a,b,m$ are integers with $m\not=0$, and $c$ and $d>0$ are rational numbers.
For this purpose, we utilize the author's philosophy for series for $1/\pi$
stated in \cite{S13} and \cite{S20}.
As a result, we make the following conjecture.

\begin{conjecture} We have
\begin{equation}\label{911}\sum_{k=0}^\infty(7k+1)\f{S_k^{(2)}(1/11)}{9^k}=\f{5445}{104\sqrt{39}\,\pi}=\f{1815\sqrt{39}}{1352\pi}
\end{equation}
and
\begin{equation}\label{18}\sum_{k=0}^\infty(1365k+181)\f{S_k^{(2)}(1/18)}{16^k}=\f{1377}{\sqrt2\,\pi}.
\end{equation}
\end{conjecture}

\begin{remark} The identities \eqref{911} and \eqref{18} are motivated by Conjectures \ref{Conj-39}
and \ref{Conj-sqrt2} respectively, and we have checked them numerically via {\tt Mathematica}.
\end{remark}

We are going to present our proofs of Theorems \ref{Th1.1}-\ref{Th1.2} and Theorem \ref{Th1.3}
in Sections 2 and 3, respectively.
In Section 4, we collect our conjectures on $p$-adic congruences involving the polynomials $S_n^{(2)}(x)\ (n\in\N)$.

\section{Proofs of Theorems \ref{Th1.1} and \ref{Th1.2}}
\setcounter{equation}{0}

As noted by Labelle \cite{Q489184},  the Chu-Vandermonde identity yields that
\begin{equation}\label{i+j}\bi{i+j}i=\sum_{k=0}^{\min\{i,j\}}\bi ik\bi jk\quad \t{for all}\ \ i,j\in\N.
\end{equation}

\noindent{\tt Proof of Theorem \ref{Th1.1}}. (i) Let $m\in\Z^+$. In light of \eqref{i+j}, we have
\begin{align*}S_n^{(m)}(x)&=\sum_{i=0}^n\sum_{j=0}^n\binom ni^m\binom nj^m\sum_{k=0}^{\min\{i,j\}}\binom ik\binom jkx^{i+j}
\\&=\sum_{k=0}^n\sum_{i=k}^n\binom ni^m\binom ikx^i\sum_{j=k}^n\binom nj^m\binom jkx^j
\\&=\sum_{k=0}^n\bigg(\sum_{i=k}^n\binom ni^{m-1}\binom nk\binom{n-k}{i-k}x^i\bigg)^2
\\&=\sum_{k=0}^n\binom nk^2x^{2k}\bigg(\sum_{j=0}^{n-k}\binom n{j+k}^{m-1}\binom{n-k}jx^j\bigg)^2.
\end{align*}
This proves \eqref{general}. Similarly, by using \eqref{i+j} we get
\begin{align*}S_n^{(0)}(x)&=\sum_{i=0}^n\sum_{j=0}^n\sum_{k=0}^{\min\{i,j\}}\bi ik\bi jkx^{i+j}
\\&=\sum_{k=0}^n\sum_{i=k}^n\bi ik x^i\sum_{j=k}^n\bi jk x^j
\\&=\sum_{k=0}^n\(\sum_{j=0}^{n-k}\bi{j+k}kx^{j+k}\)^2
\end{align*}
and hence \eqref{S0} holds.

(ii) Now we come to prove \eqref{T}. Applying \eqref{general} with $m=1$, we get
\begin{align*}S_n^{(m)}(x)&=\sum_{k=0}^n\binom nk^2x^{2k}\bigg(\sum_{j=0}^{n-k}\binom{n-k}jx^j\bigg)^2=\sum_{k=0}^n\bi nk^2x^{2k}(1+x)^{2(n-k)}.
\end{align*}

The Legendre polynomial of degree $n$ is given by
$$P_n(z)=\sum_{k=0}^n\bi nk\bi{n+k}k\l(\f{z-1}2\r)^k.$$
It is well known (cf. \cite[p.\,38, (3.134)]{G}) that
$$2^nP_n(z)=\sum_{k=0}^n\bi nk^2(z+1)^k(z-1)^{n-k}.$$
Thus
\begin{align*}&\ (2x+1)^nP_n\l(\f{2x^2+2x+1}{2x+1}\r)
\\=&\ \f{(2x+1)^n}{2^n}\sum_{k=0}^n\bi nk^2\l(\f{2(x+1)^2}{2x+1}\r)^k\l(\f{2x^2}{2x+1}\r)^{n-k}
\\=&\ \sum_{k=0}^n\bi nk^2(x+1)^{2k}(x^2)^{n-k}=\sum_{k=0}^n\bi nk^2 x^{2k}(1+x)^{2(n-k)}.
\end{align*}
Note that
\begin{align*}(2x+1)^nP_n\l(\f{2x^2+2x+1}{2x+1}\r)&=(2x+1)^n\sum_{k=0}^n\bi nk\bi{n+k}k\l(\f{x^2}{2x+1}\r)^k
\\&=\sum_{k=0}^n\bi nk\bi{n+k}kx^{2k}(2x+1)^{n-k}.
\end{align*}

Combining the last two paragraphs, we immediately obtain \eqref{T}.
Applying the Zeilberger algorithm (cf. \cite{PWZ}), we get the desired recursion formula
\eqref{rec} from \eqref{T}.
This ends our proof. \qed

For any odd prime $p$, we have
\begin{equation}\label{-1/2-k}\bi{2k}k=(-4)^k\bi{-1/2}k\eq(-4)^k\bi{(p-1)/2}k\pmod p
\end{equation}
for every $k=0,\ldots,(p-1)/2$. Thus the following lemma
follows easily from the binomial theorem.

\begin{lemma}\label{Lem-2k}
Let $p$ be an odd prime. Then
\begin{equation}\label{x^{k}}
\sum_{k=0}^{(p-1)/2}\bi{2k}kx^k\eq(1-4x)^{(p-1)/2}\pmod p
\end{equation}
and
\begin{equation}\label{kx^{k}}
\sum_{k=0}^{(p-1)/2}k\bi{2k}kx^k\eq2x(1-4x)^{(p-3)/2}\pmod p.
\end{equation}
\end{lemma}

\begin{lemma} Let $k\in\Z^+$. Then
\begin{equation}\label{ki}\f2k\sum_{k/2\ls i\ls k} i\bi ki= 2^{k-1}+\bi{2\lfloor k/2\rfloor}{\lfloor k/2\rfloor}.\end{equation}
\end{lemma}
\Proof. Observe that
\begin{align*}\f2k\sum_{k/2\ls i\ls k}i\bi ki&=2\sum_{k/2\ls i\ls k}\bi{k-1}{i-1}
=2\sum_{\lfloor(k-1)/2\rfloor\ls j\ls k-1}\bi{k-1}j
\\&=\sum_{j=\lfloor(k-1)/2\rfloor}^{k-1}\l(\bi{k-1}j+\bi{k-1}{k-1-j}\r)
\\&=\sum_{i=0}^{k-1}\bi{k-1}i+\bi{2\lfloor k/2\rfloor}{\lfloor k/2\rfloor}
=2^{k-1}+\bi{2\lfloor k/2\rfloor}{\lfloor k/2\rfloor}.
\end{align*}
This proves the desired \eqref{ki}. \qed

\medskip
\noindent{\tt Proof of Theorem \ref{Th1.2}}. Observe that
\begin{align*}\sum_{n=0}^{p-1}S_n^{(0)}(x)&=\sum_{n=0}^{p-1}\sum_{i,j=0}^n\bi{i+j}ix^{i+j}
\\&=\sum_{i,j=0}^{p-1}\bi{i+j}ix^{i+j}\sum_{n=\max\{i,j\}}^{p-1}1
=\sum_{i,j=0}^{p-1}\f{(i+j)!}{i!j!}x^{i+j}(p-\max\{i,j\})
\\&\eq-\sum_{i,j=0\atop i+j\ls p-1}^{p-1}\bi{i+j}ix^{i+j}\max\{i,j\}\pmod p.
\end{align*}
Thus
\begin{equation}\label{max}\sum_{k=0}^{p-1}S_k^{(0)}(x)\eq-\sum_{k=1}^{p-1}x^k\sum_{i=0}^k\bi ki\max\{i,k-i\}\pmod p.
\end{equation}
For any positive integer $k$, we have
\begin{align*}\sum_{i=0}^k\bi ki\max\{i,k-i\}
=&\ \sum_{k/2\ls i\ls k}\bi ki i+\sum_{k/2<j\ls k}\bi k{k-j}j
\\=&2\sum_{k/2\ls i\ls k}\bi kii-\begin{cases}\f k2\bi{k}{k/2}&\t{if}\ 2\mid k,
\\0&\t{if}\ 2\nmid k.\end{cases}
\end{align*}
Combining this with \eqref{max} and \eqref{ki}, we see that
$\sum_{k=0}^{p-1}S_k^{(0)}(x)$ is congruent to
\begin{align*}
&\ -\sum_{k=1}^{p-1}x^kk\l(2^{k-1}+\bi{2\lfloor k/2\rfloor}{\lfloor k/2\rfloor}\r)
+\sum_{k=1\atop 2\mid k}^{p-1}x^k\f k2\bi{k}{k/2}
\\=&\ -x\sum_{k=1}^{p-1}k(2x)^{k-1}-\sum_{k=1}^{p-1}kx^k\bi{2\lfloor k/2\rfloor}{\lfloor k/2\rfloor}
+\sum_{k=1}^{(p-1)/2} k\bi{2k}kx^{2k}
\end{align*}
modulo $p$.
Note that
$$\sum_{k=0}^{p}kx^k\bi{2\lfloor k/2\rfloor}{\lfloor k/2\rfloor}
=\sum_{k=0}^{(p-1)/2}\bi{2k}k\l(2kx^{2k}+(2k+1)x^{2k+1}\r).$$
Therefore, with the aid of Lemma 2.1, we have
\begin{align*}&\sum_{k=0}^{p-1}S_k^{(0)}(x)+x\sum_{k=1}^{p-1}k(2x)^{k-1}
\\\eq&\ -(2x+1)\sum_{k=0}^{(p-1)/2}k\bi{2k}kx^{2k}-x\sum_{k=0}^{(p-1)/2}\bi{2k}kx^{2k}
\\\eq&\ -(2x+1)2x^2(1-4x^2)^{(p-3)/2}-x(1-4x^2)^{(p-1)/2}
\\\eq&\begin{cases}(\f{2x^2}{2x-1}-x)(\f{1-4x^2}p)=\f x{2x-1}(\f{1-4x^2}p)\pmod p&\t{if}\ x\not\eq\f12\pmod p,\\ -\da_{p,3}\pmod p&\t{if}\ x\eq\f12\pmod p.\end{cases}
\end{align*}

If $x\eq1/2\pmod p$, then
$$\sum_{k=1}^{p-1}k(2x)^{k-1}\eq\sum_{k=1}^{p-1}k=\f{p(p-1)}2\eq0\pmod p.$$
Since
$$\sum_{k=1}^{p-1}kt^{k-1}=\f d {dt}\sum_{k=0}^{p-1}t^k=\f d{dt}\l(\f{t^p-1}{t-1}\r)
=\f{pt^{p-1}}{t-1}-\f{t^p-1}{(t-1)^2},$$
when $x\not\eq\f12\pmod p$ we obtain that
$$\sum_{k=1}^{p-1}k(2x)^{k-1}=\f{p(2x)^{p-1}}{2x-1}-\f{(2x)^p-1}{(2x-1)^2}\eq-\f1{2x-1}\pmod p$$
with the help of Fermat's little theorem. So, by the above we have the desired result. \qed

\section{Proof of Theorem \ref{Th1.3}}
\setcounter{equation}{0}

For $n\in\N$, we define the generalized central trinomial coefficient
\begin{equation}\label{Tbc}T_n(b,c)=\sum_{k=0}^{\lfloor n/2\rfloor}\bi n{2k}\bi{2k}kb^{n-2k}c^k\in\Z[b,c],
\end{equation}
which is the coefficient of $x^n$ in the expansion of $(x^2+bx+c)^n$.
It is known (cf. \cite{N}) that
$$(n+2)T_{n+2}(b,c)=(2n+3)bT_{n+1}(b,c)-(n+1)(n^2-4c)T_n(b,c)$$
for all $n\in\N$. For congruences involving $T_n(b,c)$, the reader may consult \cite{S14SCM}
The author \cite{S14,S20,S23} found nine types of series for $1/\pi$
involving generalized central trinomial coefficients.

\begin{lemma}\label{Lem2.2} For any $n\in\N$, we have
\begin{equation}\label{ST} S_n^{(1)}(x)=T_n\l(2x^2+2x+1,x^2(x+1)^2\r)\in\Z[x(x+1)].
\end{equation}
\end{lemma}
\Proof. For a polynomial $P(z)$ in $z$, let $[z^n]P(z)$ denote the coefficient of $z^n$ in the expansion of $P(z)$. Then
\begin{align*}&\ \sum_{k=0}^n\bi nk\bi n{n-k} x^{2k}((1+x)^2)^{n-k}
\\=&\ [y^n](1+x^2y)^n(1+(1+x)^2y)^{n}
\\=&\ [y^n](1+(2x^2+2x+1)y+x^2(x+1)^2y^2)^n
\\=&\ [y^n](y^{-2}+(2x^2+2x+1)y^{-1}+x^2(x+1)^2)y^{2n}
\\=&\ [z^n](z^2+(2x^2+2x+1)z+x^2(x+1)^2)
\\=&\ T_n(2x^2+2x+1,x^2(x+1)^2
\\=&\ \sum_{k=0}^{\lfloor n/2\rfloor}\bi n{2k}\bi{2k}k(2x(x+1)+1)^{n-2k}(x(x+1))^{2k}.
\end{align*}
Combining this with \eqref{T}, we immediately obtain the desired \eqref{ST}. \qed

\begin{lemma}\label{Lem2.3} Let $n\in\Z^+$. Then
\begin{equation}\label{2c}
\f{2c}n\sum_{k=0}^{n-1}T_k(b,c^2)(b-2c)^{n-1-k}
=-T_n(b,c^2)+(b+2c)T_{n-1}(b,c^2)
\end{equation}
and
\begin{equation}\label{12c^2}\begin{aligned}
&\ \f{12c^2}n\sum_{k=0}^{n-1}kT_k(b,c^2)(b-2c)^{n-1-k}-4c^2\sum_{k=0}^{n-1}T_k(b,c^2)(b-2c)^{n-1-k}
\\=&\ (b+4c)T_n(b,c^2)-(b+2c)^2T_{n-1}(b,c).
\end{aligned}
\end{equation}
\end{lemma}
\Proof. By Lemmas 3.1 and 3.2 of \cite{S14SCM}, both \eqref{2c} and \eqref{12c^2} hold for all $b,c\in\Z$. As $T_k(b,c^2)$ with $k\in\N$ is a polynomial in $b$ and $c$, we do have \eqref{2c} and \eqref{12c^2} with $b$ and $c$ as variables. This ends the proof. \qed

\medskip
\noindent{\tt Proof of Theorem \ref{Th1.3}}. (i) Let $n\in\Z^+$. In view of Lemma \ref{Lem2.2} and \eqref{Tbc},
by taking $b=2x^2+2x+1$ and $c=x(x+1)$ in \eqref{2c} we get
\begin{align*}&\ \f{2x(x+1)}n\sum_{k=0}^{n-1}S_k^{(1)}(x)
\\=&\ (2x+1)^2T_{n-1}(2x(x+1)+1,x^2(x+1)^2)-T_n(2x(x+1)+1,x^2(x+1)^2)
\end{align*}
and hence
\begin{equation}\label{S/n}\begin{aligned}&\ \f{2x(x+1)}n\sum_{k=0}^{n-1}S_k^{(1)}(x)
\\=&\ (4x(x+1)+1)\sum_{k=0}^{\lfloor (n-1)/2\rfloor}\bi {n-1}{2k}\bi{2k}k(2x(x+1)+1)^{n-1-2k}(x(x+1))^{2k}
\\&\ -\sum_{k=0}^{\lfloor n/2\rfloor}\bi n{2k}\bi{2k}k(2x(x+1)+1)^{n-2k}(x(x+1))^{2k}.
\end{aligned}\end{equation}
Thus \eqref{/n} holds provided that
\begin{equation}\label{2y}\f1{2y}\(\sum_{k=0}^{\lfloor (n-1)/2\rfloor}\bi {n-1}{2k}\bi{2k}ky^{2k}-\sum_{k=0}^{\lfloor n/2\rfloor}\bi n{2k}\bi{2k}ky^{2k}\)\in\Z[y].
\end{equation} Since $\bi{2k}k=2\bi{2k-1}{k-1}$ for all $k\in\Z^+$, we see that \eqref{2y} does hold.
So we have \eqref{/n}.

Similarly, by taking $b=2x^2+2x+1$ and $c=x(x+1)$ in \eqref{12c^2} we obtain
\begin{align*}&\ \f{12x^2(x+1)^2}n\sum_{k=0}^{n-1}kS_k^{(1)}(x)-4x^2(x+1)^2\sum_{k=0}^{n-1}S_k^{(1)}(x)
\\=&\ (6x(x+1)+1)T_n(2x(x+1)+1,x^2(x+1)^2)
\\&\ -(2x+1)^4T_{n-1}(2x(x+1)+1,x^2(x+1)^2)
\\=&\ (6x(x+1)+1)\sum_{k=0}^{\lfloor n/2\rfloor}\bi n{2k}\bi{2k}k(2x(x+1)+1)^{n-2k}(x(x+1))^{2k}
\\&\ -(2x+1)^4\sum_{k=0}^{\lfloor (n-1)/2\rfloor}\bi {n-1}{2k}\bi{2k}k(2x(x+1)+1)^{n-1-2k}(x(x+1))^{2k}.
\end{align*}
Note that
$$(2x+1)^4=(4x(x+1)+1)^2=16x^2(x+1)^2+8x(x+1)+1$$
and
$$\f{\bi{2k}k(x(x+1))^{2k}}{2x^2(x+1)^2}\in\Z[x(x+1)]\quad\t{for all}\ k\in\Z^+.$$
Also,
\begin{align*}&\ (6x(x+1)+1)(2x(x+1)+1)^n-(8x(x+1)+1)(2x(x+1)+1)^{n-1}
\\=&\ 6x(x+1)(2x(x+1)+1)^n-6x(x+1)(2x(x+1)+1)^{n-1}
\\=&\ 2x^2(x+1)^2P(x(x+1))
\end{align*}
for some polynomial $P(z)\in\Z[z]$. Thus, by the above, we see that
$$\f{6}n\sum_{k=0}^{n-1}kS_k^{(1)}(x)-2\sum_{k=0}^{n-1}S_k^{(1)}(x)\in\Z[x(x+1)].$$
Combining this with \eqref{/n} we obtain
$$\f{6}n\sum_{k=0}^{n-1}kS_k^{(1)}(x)\in\Z[x(x+1)].$$
By Lemma \ref{Lem2.2}, $S_k^{(1)}(x)\in\Z[x(x+1)]$ for all $k\in\N$.
Note also that $(6,n)=6a+nb$ for some $a,b\in\Z$. Therefore we have the desired \eqref{6/n}.

(ii) Taking $n=p$ in \eqref{S/n} we get the identity
\begin{align*}
&\ \f1p\sum_{k=0}^{p-1}S_k^{(1)}(x)
\\=&\ \f1{2x(x+1)}\sum_{k=0}^{(p-1)/2}f_p(k,x)
\bi{2k}k((2x(x+1)+1)^{p-1-2k}(x(x+1))^{2k},
\end{align*}
where
$$f_p(k,x)=(4x(x+1)+1)\bi{p-1}{2k}-(2x(x+1)+1)\bi p{2k}.$$
Note that $f_p(0,x)=4x(x+1)+1-(2x(x+1)+1)=2x(x+1)$. So
 \begin{align*}&\ \f1p\sum_{k=0}^{p-1}S_k^{(1)}(x)-(2x(x+1)+1)^{p-1}
 \\=&\ \sum_{k=1}^{(p-1)/2}f_p(k,x)
\f{\bi{2k}k}2((2x(x+1)+1)^{p-1-2k}(x(x+1))^{2k-1}.
\end{align*}
If $1\ls k\ls(p-1)/2$, then $\bi{p-1}{2k}\eq(-1)^{2k}=1\pmod p$ and $\bi p{2k}\eq0\pmod p$,
hence $f_p(k,x)\eq 4x(x+1)+1\pmod {p\Z[x(x+1)]}$. Recall that \eqref{-1/2-k} holds for each $k=1,\ldots,(p-1)/2$. So
 $$\f1p\sum_{k=0}^{p-1}S_k^{(1)}(x)-(2x(x+1)+1)^{p-1}\eq\ \f{4x(x+1)+1}2\Sigma\pmod {p\Z[x(x+1)]},$$
where
\begin{align*}\Sigma:=&\ \sum_{k=1}^{(p-1)/2}\bi{(p-1)/2}k(-4)^k((2x(x+1)+1)^{p-1-2k}(x(x+1))^{2k-1}
\\=&\ \f1{x(x+1)}\sum_{k=1}^{(p-1)/2}\bi{(p-1)/2}k\l((2x(x+1)+1)^2\r)^{(p-1)/2-k}(-4x^2(x+1)^2)^k
\\=&\ \f1{x(x+1)}\l(((2x(x+1)+1)^2-4x^2(x+1)^2)^{(p-1)/2}-(2x(x+1)+1)^{p-1}\r)
\\=&\ \f1{x(x+1)}\l(4x(x+1)+1)^{(p-1)/2}-(2x(x+1)+1)^{p-1}\r)\in 2\Z[x(x+1)].
\end{align*}
Therefore
\begin{align*}\f1p\sum_{k=0}^{p-1}S_k^{(1)}(x)
&\eq (2x(x+1)+1)^{p-1}+\f{(2x+1)^2}{2x(x+1)}\l((2x+1)^{p-1}-(2x(x+1)+1)^{p-1}\r)
\\&= \f{(2x+1)^{p+1}-(2x(x+1)+1)^p}{2x(x+1)}
\\&\eq\f{(2x+1)^{p+1}-(2x(x+1))^p-1^p}{2x(x+1)}
\\&\eq\f{(2x+1)^{p+1}-1}{2x(x+1)}-x^{p-1}(x+1)^{p-1}
\pmod{p\Z[x(x+1)]}.
\end{align*}

Observe that
\begin{align*}\f{(2x+1)^{p+1}-1}{x(x+1)}&=\l(\f1x+\f1{x+1}\r)(x+(x+1))^p+\f1{x+1}-\f1x
\\&=\l(\f1x+\f1{x+1}\r)(x^p+(x+1)^p)
\\&\quad+\l(\f1x+\f1{x+1}\r)\sum_{k=1}^{p-1}\bi pk x^k(x+1)^{p-k}+\f1{x+1}-\f1x
\\&\eq x^{p-1}+(x+1)^{p-1}+\f{x^p+1}{x+1}+\f{(x+1)^p-1}x
\\&\eq2x^{p-1}+(x+1)^{p-1}+\f{x^p+1}{x+1}\pmod{p\Z[x]}
\end{align*}
and
$$(x+1)^{p-1}=\sum_{k=0}^{p-1}\bi{p-1}kx^k\eq\sum_{k=0}^{p-1}(-x)^k=\f{x^p+1}{x+1}\pmod{p\Z[x]}.$$
Therefore
\begin{align*}\f1p\sum_{k=0}^{p-1}S_k^{(1)}(x)&\eq x^{p-1}+(x+1)^{p-1}-x^{p-1}(x+1)^{p-1}
\\&=1-(x^{p-1}-1)((x+1)^{p-1}-1)\pmod{p\Z_p[x]}.
\end{align*}
This proves \eqref{p^2}.

In view of the above, we have completed the proof of Theorem \ref{Th1.3}. \qed

 Let $p>3$ be a prime. We can also determine
$\f1p\sum_{k=0}^{p-1}kS_k^{(1)}(x)$ modulo $p\Z_p[x]$. In fact,
taking $n=p$ in the second paragraph of our proof of Theorem \ref{Th1.3}(i),
we get
\begin{align*}&\ \f{12}px^2(x+1)^2\sum_{k=0}^{p-1}kS_k^{(1)}(x)-4x^2(x+1)^2\sum_{k=0}^{p-1}S_k^{(1)}(x)
\\=&\ ((6x(x+1)+1)(2x(x+1)+1)-(2x+1)^4)(2x(x+1)+1)^{p-1}
\\&\ +\sum_{k=1}^{(p-1)/2}\bi{2k}k(2x(x+1)+1)^{p-1-2k}(x(x+1))^{2k}g_p(k,x)
\end{align*}
and hence
\begin{equation}\sigma_p:= \f 3p\sum_{k=0}^{p-1}kS_k^{(1)}(x)-\sum_{k=0}^{p-1}S_k^{(1)}(x)+(2x(x+1)+1)^{p-1}
\end{equation}
coincides with
$$\sum_{k=1}^{(p-1)/2}\f{\bi{2k-1}{k-1}}2(2x(x+1)+1)^{p-1-2k}(x(x+1))^{2(k-1)}g_p(k,x),$$
where
\begin{align*}g_p(k,x)&=(6x(x+1)+1)(2x(x+1)+1)\bi p{2k}-(2x+1)^4\bi{p-1}{2k}
\\&\eq\bi p{2k}-\bi{p-1}{2k}+(1-(2x+1)^4)\bi{p-1}{2k}
\\&\eq\bi{p-1}{2k-1}+1-(4x(x+1)+1)^2
\\&\eq p-(4x(x+1)+1)^2\pmod{2p\Z[x(x+1)]}
\end{align*}
for each $k=1,\ldots,(p-1)/2$. (Note that $\bi{p-1}{2k-1}=\f{p-1}{2k-1}\bi{p-2}{2k-2}$
is an even number congruent to $(-1)^{2k-1}=-1$ modulo $p$.)
It follows that $\sigma_p$ is congruent to
$$(p-(2x+1)^4)\sum_{k=1}^{(p-1)/2}\f{\bi{2k}k}4(2x(x+1)+1)^{p-1-2k}(x(x+1))^{2(k-1)}$$
modulo $p\Z[x(x+1)]$.
Recall that \eqref{-1/2-k} holds for all $k=1,\ldots,(p-1)/2$. Thus $\sigma_p$ is congruent to
\begin{align*}&\ -\f{(2x+1)^4}{4x^2(x+1)^2}\sum_{k=1}^{(p-1)/2}
\bi{(p-1)/2}k\l((2x(x+1)+1)^2\r)^{(p-1)/2-k}(-4x^2(x+1)^2)^k
\\=&\ -\f{(2x+1)^4}{4x^2(x+1)^2}\l(((2x(x+1)+1)^2-4x^2(x+1)^2)^{(p-1)/2}-(2x(x+1)+1)^{p-1}\r)
\\=&\ -\f{(2x+1)^4}{4x^2(x+1)^2}\l((2x+1)^{p-1}-(2x(x+1)+1)^{p-1}\r)
\end{align*}
modulo $p\Z_p[x(x+1)]$.
Combining this with \eqref{p^2}, we get
\begin{equation} \f 3p\sum_{k=0}^{p-1}kS_k^{(1)}(x)
\eq -\f{(2x+1)^{p+3}}{4x^2(x+1)^2}+\f{6x(x+1)+1}{4x^2(x+1)^2}(2x(x+1)+1)^p\pmod{p\Z_p[x]}.
\end{equation}

\section{Conjectures on congruences involving $S_n^{(2)}(x)$}
 \setcounter{equation}{0}

\begin{conjecture} Let $p$ be an odd prime.

{\rm (i)} We have
$$\sum_{k=0}^{p-1}2^kS_k^{(2)}\l(-\f12\r)\eq
\begin{cases}4x^2-2p\pmod{p^2}&\t{if}\ p=x^2+y^2\ (x,y\in\Z)\ \t{with}\ 2\nmid x,
\\0\pmod{p^2}&\t{if}\ p\eq3\pmod4.\end{cases}$$

{\rm (ii)} We have
$$\sum_{k=0}^{p-1}(3k+2)2^kS_k^{(2)}\l(-\f12\r)\eq 2p\l(2-\l(\f 2p\r)\r)\pmod{p^2}.$$
Moreover, when $p\eq\pm1\pmod8$ we have
$$\f1{(pn)^2}\bigg(\sum_{k=0}^{pn-1}(3k+2)2^kS_k^{(2)}\l(-\f12\r)-p\sum_{k=0}^{n-1}(3k+2)2^k
S_k^{(2)}\l(-\f12\r)\bigg)
\in\Z_p.$$
\end{conjecture}
\begin{remark} By elementary number theory, each prime $p\eq1\pmod4$
can be written uniquely as $x^2+y^2$ with $x,y\in\Z^+$, $2\nmid x$ and $2\mid y$.
\end{remark}

\begin{conjecture} Let $p>3$ be a prime.

{\rm (i)} We have
$$\sum_{k=0}^{p-1}\frac{S_k^{(2)}(-1)}{(-3)^k}
\equiv\begin{cases}4x^2-2p\pmod{p^2}&\text{if}\ p=x^2+y^2\ \text{with}\ 3\nmid x\ \text{and}\ 3\mid y,\\4xy\pmod{p^2}&\text{if}\ p=x^2+y^2\ \text{with}\ x\equiv y\not\equiv0\pmod3,\\0\pmod{p^2}&\text{if}\ p\equiv3\pmod4,\end{cases}$$
where $x$ and $y$ are integers.

{\rm (ii)} We have the congruence
$$\sum_{k=0}^{p-1}(28k+17)\frac{S_k^{(2)}(-1)}{(-3)^k}\equiv p\left(11+6\left(\frac p3\right)\right)\pmod{p^2}.$$
Moreover, when $p\equiv1\pmod3$, for any positive integer $n$ we have
$$\frac1{(pn)^2}\bigg(\sum_{k=0}^{pn-1}(28k+17)\frac{S_k^{(2)}(-1)}{(-3)^k}
-p\sum_{k=0}^{n-1}(28k+17)\frac{S_k^{(2)}(-1)}{(-3)^k}\bigg)\in\mathbb Z_p.$$
\end{conjecture}
\begin{remark} $S_0^{(2)}(-1),\ldots,S_9^{(2)}(-1)$ take the values
1, 1, 9, 73, 361, 5001, 35001, 348489, 3693033, 31360681, respectively.
\end{remark}

\begin{conjecture} Let $p$ be an odd prime.

{\rm (i)} We have
\begin{align*}&\sum_{k=0}^{p-1}9^kS_k^{(2)}\l(\f23\r)\eq\sum_{k=0}^{p-1}9^kS_k^{(2)}\l(\f43\r)
\\\equiv&\begin{cases}4x^2-2p\pmod{p^2}&\text{if}\ p=x^2+y^2\ \text{with}\ 3\nmid x\ \text{and}\ 3\mid y,\\4xy\pmod{p^2}&\text{if}\ p=x^2+y^2\ \text{with}\ x\equiv y\not\equiv0\pmod3,\\0\pmod{p^2}&\text{if}\ p\equiv3\pmod4,\end{cases}
\end{align*}
where $x$ and $y$ are integers.

{\rm (ii)} For any $n\in\Z^+$ we have
$$\f2n\sum_{k=0}^{n-1}(3k+2)9^kS_k^{(2)}\l(\f 23\r)\in\Z$$
and
$$\frac1{(pn)^2}\bigg(\sum_{k=0}^{pn-1}(3k+2)9^kS_k^{(2)}\l(\f23\r)
-p\sum_{k=0}^{n-1}(3k+2)9^kS_k^{(2)}\l(\f23\r)\bigg)\in\mathbb Z_p.$$
When $p\eq1\pmod4$, we also have
$$\sum_{k=0}^{p-1}(24k+19)9^kS_k^{(2)}\l(\f 43\r)\eq p\pmod{p^2}.$$
\end{conjecture}
\begin{remark} For $a\in\Z$ and $n\in\N$, we clearly have
$$9^nS_n^{(2)}\l(\f a3\r)=\sum_{i,j=0}^n\bi ni^2\bi nj^2 a^{i+j}3^{2n-i-j}\in\Z.$$
\end{remark}

\begin{conjecture}\label{Conj4.4} Let $p>5$ be a prime.

{\rm (i)} We have
\begin{align*}&\sum_{k=0}^{p-1}25^kS_k^{(2)}\l(\f4{15}\r)
\\\equiv&\begin{cases}4x^2-2p\pmod{p^2}&\text{if}\ p=x^2+y^2\ (x,y\in\Z,\ 5\nmid x\ \&\ 5\mid y),
\\4xy\pmod{p^2}&\text{if}\  p=x^2+y^2\ (x,y\in\Z\ \&\ x\eq y\not\eq0\pmod5),
\\0\pmod{p^2}&\text{if}\ p\eq3\pmod4.
\end{cases}\end{align*}

{\rm (ii)} For any $n\in\Z^+$ we have
$$\f1{(pn)^2}\(\sum_{k=0}^{pn-1}(168k+125)25^kS_k^{(2)}\l(\f4{15}\r)-p\sum_{k=0}^{n-1}
(168k+125)25^kS_k^{(2)}\l(\f4{15}\r)\)\in\Z_p.$$
\end{conjecture}
\begin{remark} For any prime $p=x^2+y^2>5$ with $x,y\in\Z$ and $5\nmid xy$, either $x^2\eq y^2\eq1\pmod 5$ or $x^2\eq y^2\eq 4\pmod5$, and hence $x$ is congruent to $y$ or $-y$ modulo $5$.
\end{remark}

\begin{conjecture} Let $p>7$ be a prime.

{\rm (i)} We have
\begin{align*}&\sum_{k=0}^{p-1}\f{S_k^{(2)}(1/7)}{5^k}
\\\equiv&\begin{cases}4x^2-2p\pmod{p^2}&\text{if}\ p=x^2+y^2\ (x,y\in\Z,\ 5\nmid x\ \&\ 5\mid y),
\\4xy\pmod{p^2}&\text{if}\  p=x^2+y^2\ (x,y\in\Z\ \&\ x\eq y\not\eq0\pmod5),
\\0\pmod{p^2}&\text{if}\ p\eq3\pmod4.
\end{cases}\end{align*}

{\rm (ii)} We have the congruence
$$\sum_{k=0}^{p-1}(124k+43)\f{S_k^{(2)}(1/7)}{5^k}\eq p\l(53-10\l(\f p5\r)\r)\pmod{p^2}.$$
Moreover, when $p\eq1,4\pmod5$, for any $n\in\Z^+$ we have
$$\f1{(pn)^2}\(\sum_{k=0}^{pn-1}(124k+43)\f{S_k^{(2)}(1/7)}{5^k}-p\sum_{k=0}^{n-1}
(124k+43)\f{S_k^{(2)}(1/7)}{5^k}\)\in\Z_p.$$
\end{conjecture}
\begin{remark} This is similar to Conjecture \ref{Conj4.4}.
\end{remark}

\begin{conjecture} Let $p$ be an odd prime.

{\rm (i)} We have
$$\begin{aligned}&\sum_{k=0}^{p-1}16^kS_k^{(2)}\l(\f{3}{8}\r)\eq\sum_{k=0}^{p-1}81^kS_k^{(2)}\l(\f 89\r)
\\\equiv&\ \begin{cases}(\f p3)(4x^2-2p)\pmod{p^2}&\text{if}\  p=x^2+2y^2\ \t{with}\  x,y\in\Z,\\0\pmod{p^2}&\text{if}\ p\equiv 5,7\pmod{8}.\end{cases}\end{aligned}$$
If $p\not=5$, then
$$\begin{aligned}\sum_{k=0}^{p-1}S_k^{(2)}\l(\f{12}{5}\r)
\equiv&\ \begin{cases}(\f p3)(4x^2-2p)\pmod{p^2}&\text{if}\  p=x^2+2y^2\ \t{with}\  x,y\in\Z,\\0\pmod{p^2}&\text{if}\ p\equiv 5,7\pmod{8}.\end{cases}\end{aligned}$$

{\rm (ii)} When $p\eq1,3\pmod8$, we have
$$\sum_{k=0}^{p-1}(480k+377)81^kS_k^{(2)}\l(\f 89\r)\eq 122p\pmod{p^2}.$$
If $p>3$ and $n\in\Z^+$, then
$$\f1{(pn)^2}\bigg(\sum_{k=0}^{pn-1}(32k+23)16^kS_k^{(2)}\l(\f{3}{8}\r)-p\sum_{k=0}^{n-1}(32k+23)16^kS_k^{(2)}
\l(\f{3}{8}\r)\bigg)\in\Z_p.$$
If $p\not=5$ and $n\in\Z^+$, then
$$\f1{(pn)^2}\bigg(\sum_{k=0}^{pn-1}(104k+69)S_k^{(2)}\l(\f{12}{5}\r)-p\sum_{k=0}^{n-1}(104k+69)S_k^{(2)}
\l(\f{12}{5}\r)\bigg)\in\Z_p.$$

\end{conjecture}
\begin{remark} It is well known (cf. \cite{Cox}) that any prime $p\eq1,3\pmod 8$ can be written uniquely as $x^2+2y^2$ with $x,y\in\Z^+$.
\end{remark}

\begin{conjecture} \label{Conj4.7} Let $p$ be an odd prime.

{\rm (i)} We have
$$\begin{aligned}\sum_{k=0}^{p-1}S_k^{(2)}(2)
\equiv&\ \begin{cases}4x^2-2p\pmod{p^2}&\text{if}\ p\equiv1,9\pmod{20}\ \&\ p=x^2+5y^2,\\2x^2-2p\pmod{p^2}&\text{if}\ p\equiv3,7\pmod{20}\ \&\ 2p=x^2+5y^2,\\0\pmod{p^2}&\text{if}\ p\equiv 11,13,17,19\pmod{20},\end{cases}\end{aligned}$$
where $x$ and $y$ are integers.
Also,
\begin{align*}&\sum_{k=0}^{p-1}81^kS_k^{(2)}\l(\f{10}9\r)
\\\equiv&\ \begin{cases}(\f p3)(4x^2-2p)\pmod{p^2}&\text{if}\ p\equiv1,9\pmod{20}\ \&\ p=x^2+5y^2,\\(\f p3)(2p-2x^2)\pmod{p^2}&\text{if}\ p\equiv3,7\pmod{20}\ \&\ 2p=x^2+5y^2,\\0\pmod{p^2}&\text{if}\ p\equiv 11,13,17,19\pmod{20},\end{cases}\end{align*}
where $x$ and $y$ are integers.

{\rm (ii)} If $(\frac{-5}p)=1$ $($i.e., $p\equiv 1,3,7,9\pmod{20})$, then
$$\sum_{k=0}^{p-1}(8k+5)S_k^{(2)}(2)\equiv\f{14}5p\pmod {p^2}.$$
Also, we have
$$\sum_{k=0}^{p-1}(15k+13)81^kS_k^{(2)}\l(\f{10}9\r)\eq\f p{16}\l(43+165\l(\f{-15}p\r)\r)\pmod{p^2}.$$
Moreover, when $(\f{-15}p)=1$, for any $n\in\Z^+$ we have
$$\f1{(pn)^2}\bigg(\sum_{k=0}^{pn-1}(15k+13)81^kS_k^{(2)}\l(\f{10}9\r)
-p\sum_{k=0}^{n-1}(15k+13)81^kS_k^{(2)}\l(\f{10}9\r)\bigg)\in\Z_p.$$
\end{conjecture}
\begin{remark} It is known (cf. \cite{Cox}) that any prime $p\eq1,9\pmod{20}$ can be written uniquely
as $x^2+5y^2$ with $x,y\in\Z^+$, and for any prime $p\eq3,7\pmod{20}$ we can write $2p$
as $x^2+5y^2$ with $x,y\in\Z^+$ in a unique way.
\end{remark}

\begin{conjecture} Let $p>3$ be a prime with $p\not=11$.

{\rm (i)} We have
$$\begin{aligned}\sum_{k=0}^{p-1}S_k^{(2)}\l(-\f2{11}\r)
\equiv& \begin{cases}4x^2-2p\pmod{p^2}&\text{if}\ p\equiv1,9\pmod{20}\ \&\ p=x^2+5y^2,\\2x^2-2p\pmod{p^2}&\text{if}\ p\equiv3,7\pmod{20}\ \&\ 2p=x^2+5y^2,\\0\pmod{p^2}&\text{if}\ p\equiv 11,13,17,19\pmod{20},\end{cases}\end{aligned}$$
where $x$ and $y$ are integers.

{\rm (ii)} We have the congruence
$$\sum_{k=0}^{p-1}(100k+61)S_k^{(2)}\l(-\f2{11}\r)\eq\f p{11}\l(826\l(\f 3p\r)-155\r)\pmod{p^2}.$$
Moreover, when $(\frac{3}p)=1$ $($i.e., $p\equiv \pm1\pmod{12})$, for any $n\in\Z^+$ we have
$$\f1{(pn)^2}\bigg(\sum_{k=0}^{pn-1}(100k+61)S_k^{(2)}\l(-\f2{11}\r)-p\sum_{k=0}^{n-1}(100k+61)S_k^{(2)}
\l(-\f2{11}\r)\bigg)\in\Z_p.$$
\end{conjecture}
\begin{remark} This is similar to Conjecture \ref{Conj4.7}.
\end{remark}

\begin{conjecture} Let $p>3$ be a prime.

{\rm (i)} We have
\begin{align*}&\sum_{k=0}^{p-1}S_k^{(2)}\l(\f13\r)\eq\sum_{k=0}^{p-1}S_k^{(2)}\l(\f43\r)
\\\equiv&\begin{cases}4x^2-2p\pmod{p^2}&\text{if}\ (\frac{-2}p)=(\frac p{5})=1
\ \&\ p=x^2+10y^2,
\\8x^2-2p\pmod{p^2}&\text{if}\ (\frac{-2}{p})=(\frac p{5})=-1\ \&\ p=2x^2+5y^2,
\\0\pmod{p^2}&\text{if}\ (\frac{-10}p)=-1,
\end{cases}\end{align*}
where $x$ and $y$ are integers.

{\rm (ii)} We have the congruence
$$\sum_{k=0}^{p-1}(60k+31)S_k^{(2)}\l(\f 13\r)\eq p\l(41-10\l(\f p5\r)\r)\pmod{p^2}.$$
Moreover, when $p\eq1,4\pmod5$, for any $n\in\Z^+$ we have
$$\f1{(pn)^2}\(\sum_{k=0}^{pn-1}(60k+31)S_k^{(2)}\l(\f13\r)-p\sum_{k=0}^{n-1}(60k+31)
S_k^{(2)}\l(\f13\r)\)\in\Z_p.$$
If $(\f{-10}p)=1$, then
$$\sum_{k=0}^{p-1}(25k+13)S_k^{(2)}\l(\f 43\r)\eq\f{245}{24}p\pmod{p^2}.$$
\end{conjecture}
\begin{remark} The imaginary quadratic field $\mathbb Q(\sqrt{-10})$ has class number $2$.
\end{remark}

\begin{conjecture} Let $p>3$ be a prime.

{\rm (i)} We have
\begin{align*}&\sum_{k=0}^{p-1}S_k^{(2)}\l(-\f23\r)
\eq\sum_{k=0}^{p-1}9^kS_k^{(2)}\l(-\f23\r)\eq\sum_{k=0}^{p-1}\frac{S_k^{(2)}(-1/2)}{(-4)^k}
\eq\sum_{k=0}^{p-1}16^kS_k^{(2)}\l(\f 98\r)
\\\equiv&\begin{cases}4x^2-2p\pmod{p^2}&\text{if}\ (\frac{-1}p)=(\frac p{13})=1
\ \&\ p=x^2+13y^2,
\\2x^2-2p\pmod{p^2}&\text{if}\ (\frac{-1}{p})=(\frac p{13})=-1\ \&\ 2p=x^2+13y^2,
\\0\pmod{p^2}&\text{if}\ (\frac{-13}p)=-1,
\end{cases}\end{align*}
where $x$ and $y$ are integers.

{\rm (ii)} Let $n\in\Z^+$. Then
$$\f1{(pn)^2}\bigg(\sum_{k=0}^{pn-1}(39k+25)S_k^{(2)}\l(-\f23\r)-p\sum_{k=0}^{n-1}(39k+25)S_k^{(2)}\l(-\f23\r)\bigg)\in\Z_p$$
and
$$\f1n\sum_{k=0}^{n-1}(6k+5)9^kS_k^{(2)}\l(-\f23\r)\in\Z^+.$$
When $(\f{-13}p)=1$, we also have
$$\f1{n^2}\sum_{k=0}^{n-1}(1872k+1387)16^kS_k^{(2)}\l(\f98\r)\in\Z_p.$$
\end{conjecture}
\begin{remark} The imaginary quadratic field $\mathbb Q(\sqrt{-13})$ has class number $2$.
\end{remark}

\begin{conjecture} Let $p>3$ be a prime.

{\rm (i)} We have
\begin{align*}&\sum_{k=0}^{p-1}\frac{S_k^{(2)}(2)}{9^k}
\\\equiv&\begin{cases}4x^2-2p\pmod{p^2}&\text{if}\ (\frac{-1}p)=(\frac p3)=(\frac p7)=1
\ \&\ p=x^2+21y^2,
\\2x^2-2p\pmod{p^2}&\text{if}\ (\frac p7)=1,\ (\frac{-1}p)=(\frac p3)=-1\ \&\ 2p=x^2+21y^2,
\\12x^2-2p\pmod{p^2}&\text{if}\ (\frac p3)=1,\ (\frac{-1}p)=(\frac p7)=-1\ \&\ p =3x^2+7y^2,
\\6x^2-2p\pmod{p^2}&\text{if}\ (\frac{-1}p)=1,\ (\frac p3)=(\frac p7)=-1\ \&\ 2p=3x^2+7y^2,
\\0\pmod{p^2}&\text{if}\ (\frac{-21}p)=-1,
\end{cases}\end{align*}
where $x$ and $y$ are integers.

{\rm (ii)} For any positive integer $n$, we have
$$\frac1{(pn)^2}\bigg(\sum_{k=0}^{pn-1}(7k+3)\frac{S_k^{(2)}(2)}{9^k}-p\sum_{k=0}^{n-1}(7k+3)
\frac{S_k^{(2)}(2)}{9^k}\bigg)\in\mathbb Z_p.$$
\end{conjecture}
\begin{remark} The imaginary quadratic field $\mathbb Q(\sqrt{-21})$ has class number $4$.
\end{remark}

\begin{conjecture} \label{Conj4.12} Let $p>3$ be a prime.

{\rm (i)} We have
\begin{align*}\sum_{k=0}^{p-1}\frac{S_k^{(2)}(-4)}{81^k}\eq&\sum_{k=0}^{p-1}\f{S_k^{(2)}(1/6)}{4^k}
\\\equiv&\begin{cases}4x^2-2p\pmod{p^2}&\text{if}\ (\frac{2}p)=(\frac p{11})=1
\ \&\ p=x^2+22y^2,
\\8x^2-2p\pmod{p^2}&\text{if}\ (\frac2{p})=(\frac p{11})=-1\ \&\ p=2x^2+11y^2,
\\0\pmod{p^2}&\text{if}\ (\frac{-22}p)=-1,
\end{cases}\end{align*}
where $x$ and $y$ are integers.

{\rm (ii)} For any positive integer $n$, we have
$$\frac1{(pn)^2}\bigg(\sum_{k=0}^{pn-1}(11k+6)\frac{S_k^{(2)}(-4)}{81^k}-p
\sum_{k=0}^{n-1}(11k+6)\frac{S_k^{(2)}(-4)}{81^k}\bigg)\in\mathbb Z_p.$$
\end{conjecture}
\begin{remark} The imaginary quadratic field $\mathbb Q(\sqrt{-22})$ has class number $2$.
\end{remark}

\begin{conjecture} Let $p>5$ be a prime.

{\rm (i)} We have
\begin{align*}&\sum_{k=0}^{p-1}S_k^{(2)}\l(\f8{15}\r)
\equiv\begin{cases}4x^2-2p\pmod{p^2}&\text{if}\ (\frac{2}p)=(\frac p{11})=1
\ \&\ p=x^2+22y^2,
\\8x^2-2p\pmod{p^2}&\text{if}\ (\frac2{p})=(\frac p{11})=-1\ \&\ p=2x^2+11y^2,
\\0\pmod{p^2}&\text{if}\ (\frac{-22}p)=-1,
\end{cases}\end{align*}
where $x$ and $y$ are integers.

{\rm (ii)} For any positive integer $n$, we have
$$\frac1{(pn)^2}\bigg(\sum_{k=0}^{pn-1}(2244k+1147)S_k^{(2)}\l(\f8{15}\r)-p
\sum_{k=0}^{n-1}(2244k+1147)S_k^{(2)}\l(\f8{15}\r)\bigg)\in\mathbb Z_p.$$
\end{conjecture}
\begin{remark} This is similar to Conjecture \ref{Conj4.12}.
\end{remark}

\begin{conjecture}\label{Conj4.14} Let $p>3$ be a prime.

{\rm (i)} We have
\begin{align*}&\sum_{k=0}^{p-1}S_k^{(2)}\l(-\f 43\r)
\\\eq&\begin{cases}4x^2-2p\pmod{p^2}&\t{if}\ (\f 2p)=(\f p3)-(\f p5)=1\ \&\ p=x^2+30y^2,
\\8x^2-2p\pmod{p^2}&\t{if}\ (\f 2p)=1,\ (\f p3)=(\f p5)=-1\ \&\ p=2x^2+15y^2,
\\12x^2-2p\pmod{p^2}&\t{if}\ (\f p3)=1,\ (\f 2p)=(\f p5)=-1\ \&\ p=3x^2+10y^2,
\\2p-20x^2\pmod{p^2}&\t{if}\ (\f p5)=1,\ (\f 2p)=(\f p3)=-1\ \&\ p=5x^2+6y^2,
\\0\pmod{p^2}&\t{if}\ (\f{-30}p)=-1,
\end{cases}\end{align*}
where $x$ and $y$ are integers.

{\rm (ii)} We have the congruence
$$\sum_{k=0}^{p-1}(120k+79)S_k^{(2)}\l(-\f 43\r)
\eq p\l(47+32\l(\f 2p\r)\r)\pmod{p^2}.$$
Moreover, when $p\eq\pm1\pmod8$, for any $n\in\Z^+$ we have
$$\f1{(pn)^2}\bigg(\sum_{k=0}^{pn-1}(120k+79)S_k^{(2)}\l(-\f 43\r)
-p\sum_{k=0}^{n-1}(120k+79)S_k^{(2)}\l(-\f 43\r)\bigg)\in\Z_p.$$
\end{conjecture}
\begin{remark} The imaginary quadratic field $\Q(\sqrt{-30})$ has class number $4$.
\end{remark}

\begin{conjecture} Let $p>3$ be a prime with $p\not=11$.

{\rm (i)} We have
\begin{align*}&\sum_{k=0}^{p-1}\frac{S_k^{(2)}(-4)}{9^k}
\\\equiv&\begin{cases}4x^2-2p\pmod{p^2}&\text{if}\ (\frac{-1}p)=(\frac p3)=(\frac p{11})=1
\ \&\ p=x^2+33y^2,
\\2x^2-2p\pmod{p^2}&\text{if}\ (\frac {-1}p)=1,\ (\frac{p}3)=(\frac p{11})=-1\ \&\ 2p=x^2+33y^2,
\\12x^2-2p\pmod{p^2}&\text{if}\ (\frac p{11})=1,\ (\frac{-1}p)=(\frac p3)=-1\ \&\ p =3x^2+11y^2,
\\6x^2-2p\pmod{p^2}&\text{if}\ (\frac p3)=1,\ (\frac {-1}p)=(\frac p{11})=-1\ \&\ 2p=3x^2+11y^2,
\\0\pmod{p^2}&\text{if}\ (\frac{-33}p)=-1,
\end{cases}\end{align*}
where $x$ and $y$ are integers.

{\rm (ii)} For any positive integer $n$, we have
$$\frac1{(pn)^2}\bigg(\sum_{k=0}^{pn-1}(8k+5)\frac{S_k^{(2)}(-4)}{9^k}-p\left(\frac{33}p\right)
\sum_{k=0}^{n-1}(8k+5)\frac{S_k^{(2)}(-4)}{9^k}\bigg)\in\mathbb Z_p.$$
\end{conjecture}
\begin{remark} The imaginary quadratic field $\mathbb Q(\sqrt{-33})$ has class number $4$.
\end{remark}

\begin{conjecture} Let $p>5$ be a prime.

{\rm (i)} We have
\begin{align*}&\sum_{k=0}^{p-1}\frac{S_k^{(2)}(18)}{625^k}\eq\sum_{k=0}^{p-1}\frac{S_k^{(2)}(16/5)}{81^k}
\\\equiv&\begin{cases}4x^2-2p\pmod{p^2}&\text{if}\ (\frac{-1}p)=(\frac p{37})=1
\ \&\ p=x^2+37y^2,
\\2x^2-2p\pmod{p^2}&\text{if}\ (\frac{-1}{p})=(\frac p{37})=-1\ \&\ 2p=x^2+37y^2,
\\0\pmod{p^2}&\text{if}\ (\frac{-37}p)=-1,
\end{cases}\end{align*}
where $x$ and $y$ are integers.

{\rm (ii)} If $(\f{-37}p)=1$, then
$$\sum_{k=0}^{p-1}(263736k+173561)\f{S_k^{(2)}(18)}{625^k}\eq \f{386918}5p\pmod {p^2}.$$
For any positive integer $n$, we have
$$\frac1{(pn)^2}\bigg(\sum_{k=0}^{pn-1}(53872k+20157)\frac{S_k^{(2)}(16/5)}{81^k}-p
\sum_{k=0}^{n-1}(53872k+20157)\frac{S_k^{(2)}(16/5)}{81^k}\bigg)\in\mathbb Z_p.$$
\end{conjecture}
\begin{remark} The imaginary quadratic field $\mathbb Q(\sqrt{-37})$ has class number $2$.
\end{remark}

\begin{conjecture} Let $p$ be an odd prime.

{\rm (i)} We have
\begin{align*}&\sum_{k=0}^{p-1}81^kS_k^{(2)}\l(\f 49\r)
\\\eq&\begin{cases}4x^2-2p\pmod{p^2}&\t{if}\ (\f {-2}p)=(\f p3)=(\f p{7})=1\ \&\ p=x^2+42y^2,
\\ 8x^2-2p\pmod{p^2}&\t{if}\ (\f {p}7)=1,\ (\f {-2}p)=(\f p{3})=-1\ \&\ p=2x^2+21y^2,
\\12x^2-2p\pmod{p^2}&\t{if}\ (\f {-2}p)=1,\ (\f {p}3)=(\f p{7})=-1\ \&\ p=3x^2+14y^2,
\\24x^2-2p\pmod{p^2}&\t{if}\ (\f p{3})=1,\ (\f {-2}p)=(\f p7)=-1\ \&\ p=6x^2+7y^2,
\\0\pmod{p^2}&\t{if}\ (\f{-42}p)=-1,
\end{cases}\end{align*}
where $x$ and $y$ are integers.

{\rm (ii)} We have the congruence
$$\sum_{k=0}^{p-1}(252k+185)81^kS_k^{(2)}\l(\f49\r)
\eq p\l(17+168\l(\f {6}p\r)\r)\pmod{p^2}.$$
Moreover, when $(\f{6}p)=1$, for any $n\in\Z^+$ we have
$$\f1{(pn)^2}\bigg(\sum_{k=0}^{pn-1}(252k+185)81^kS_k^{(2)}\l(\f 49\r)
-p\sum_{k=0}^{n-1}(252k+185)81^kS_k^{(2)}\l(\f 49\r)\bigg)\in\Z_p.$$
\end{conjecture}
\begin{remark} The imaginary quadratic field $\Q(\sqrt{-42})$ has class number $4$.
\end{remark}

\begin{conjecture} Let $p>3$ be a prime with $p\not=7$.

{\rm (i)} We have
\begin{align*}&\sum_{k=0}^{p-1}\f{S_k^{(2)}(2/7)}{9^k}
\\\eq&\begin{cases}4x^2-2p\pmod{p^2}&\t{if}\ (\f {-1}p)=(\f p3)=(\f p{19})=1\ \&\ p=x^2+57y^2,
\\ 2x^2-2p\pmod{p^2}&\t{if}\ (\f {-1}p)=1,\ (\f p5)=(\f p{19})=-1\ \&\ 2p=x^2+57y^2,
\\12x^2-2p\pmod{p^2}&\t{if}\ (\f p3)=1,\ (\f {-1}p)=(\f p{19})=-1\ \&\ p=3x^2+19y^2,
\\6x^2-2p\pmod{p^2}&\t{if}\ (\f p{19})=1,\ (\f {-1}p)=(\f p3)=-1\ \&\ 2p=3x^2+19y^2,
\\0\pmod{p^2}&\t{if}\ (\f{-57}p)=-1,
\end{cases}\end{align*}
where $x$ and $y$ are integers.

{\rm (ii)} We have the congruence
$$\sum_{k=0}^{p-1}(133k+48)\f{S_k^{(2)}(2/7)}{9^k}
\eq \f p{14}\l(1185-513\l(\f {57}p\r)\r)\pmod{p^2}.$$
Moreover, when $(\f{57}p)=1$, for any $n\in\Z^+$ we have
$$\f1{(pn)^2}\bigg(\sum_{k=0}^{pn-1}(133k+48)\f{S_k^{(2)}(2/7)}{9^k}
-p\sum_{k=0}^{n-1}(133k+48)\f{S_k^{(2)}(2/7)}{9^k}\bigg)\in\Z_p.$$
\end{conjecture}
\begin{remark} The imaginary quadratic field $\Q(\sqrt{-57})$ has class number $4$.
\end{remark}

\begin{conjecture} Let $p>3$ be a prime.

{\rm (i)} We have
\begin{align*}&\sum_{k=0}^{p-1}\l(\f{16}{81}\r)^{k}S_k^{(2)}\l(-\f 58\r)
\\\eq&\begin{cases}4x^2-2p\pmod{p^2}&\t{if}\ (\f 2p)=(\f p5)=(\f p7)=1\ \&\ p=x^2+70y^2,
\\ 8x^2-2p\pmod{p^2}&\t{if}\ (\f p7)=1,\ (\f 2p)=(\f p5)=-1\ \&\ p=2x^2+35y^2,
\\20x^2-2p\pmod{p^2}&\t{if}\ (\f p5)=1,\ (\f 2p)=(\f p7)=-1\ \&\ p=5x^2+14y^2,
\\2p-28x^2\pmod{p^2}&\t{if}\ (\f 2p)=1,\ (\f p5)=(\f p7)=-1\ \&\ p=7x^2+10y^2,
\\0\pmod{p^2}&\t{if}\ (\f{-70}p)=-1,
\end{cases}\end{align*}
where $x$ and $y$ are integers.

{\rm (ii)} We have the congruence
$$\sum_{k=0}^{p-1}(560k+303)\l(\f{16}{81}\r)^kS_k^{(2)}\l(-\f 58\r)
\eq 303p\pmod{p^2}.$$
\end{conjecture}
\begin{remark} The imaginary quadratic field $\Q(\sqrt{-70})$ has class number $4$.
\end{remark}

\begin{conjecture} \label{Conj-39} Let $p>3$ be a prime with $p\not=11$.

{\rm (i)} We have
\begin{align*}&\sum_{k=0}^{p-1}\f{S_k^{(2)}(1/11)}{9^k}
\\\equiv&\begin{cases}4x^2-2p\pmod{p^2}&\text{if}\ (\f 2p)=(\f p3)=(\f p{13})=1\ \&\ p=x^2+78y^2,
\\8x^2-2p\pmod{p^2}&\text{if}\  (\f 2p)=1,\ (\f p3)=(\f p{13})=-1\ \&\ p=2x^2+39y^2,
\\12x^2-2p\pmod{p^2}&\text{if}\  (\f p{13})=1,\ (\f 2p)=(\f p{3})=-1\ \&\ p=3x^2+26y^2,
\\24x^2-2p\pmod{p^2}&\text{if}\  (\f p{3})=1,\ (\f 2p)=(\f p{13})=-1\ \&\ p=6x^2+13y^2,
\\0\pmod{p^2}&\text{if}\ (\f{-78}p)=-1,
\end{cases}\end{align*}
where $x$ and $y$ are integers.

{\rm (ii)} We have the congruence
$$\sum_{k=0}^{p-1}(7k+1)\f{S_k^{(2)}(1/11)}{9^k}\eq \f p{26}\l(35\l(\f{-39}p\r)-9\l(\f {13}p\r)\r)\pmod{p^2}.$$
Moreover, when $p\eq1\pmod3$, for any $n\in\Z^+$ we have
$$\f1{(pn)^2}\(\sum_{k=0}^{pn-1}(7k+1)\f{S_k^{(2)}(1/11)}{9^k}-p\l(\f p{13}\r)\sum_{k=0}^{n-1}
(7k+1)\f{S_k^{(2)}(1/11)}{9^k}\)\in\Z_p.$$
\end{conjecture}
\begin{remark} The imaginary quadratic field $\mathbb Q(\sqrt{-78})$ has class number $4$.
\end{remark}

\begin{conjecture} \label{Conj-sqrt2} Let $p>3$ be a prime.

{\rm (i)} We have
\begin{align*}&\sum_{k=0}^{p-1}\f{S_k^{(2)}(1/18)}{16^k}
\\\equiv&\begin{cases}4x^2-2p\pmod{p^2}&\text{if}\ (\f {-2}p)=(\f p5)=(\f p{13})=1\ \&\ p=x^2+130y^2,
\\8x^2-2p\pmod{p^2}&\text{if}\  (\f {-2}p)=1,\ (\f p5)=(\f p{13})=-1\ \&\ p=2x^2+65y^2,
\\20x^2-2p\pmod{p^2}&\text{if}\  (\f p{5})=1,\ (\f {-2}p)=(\f p{13})=-1\ \&\ p=5x^2+26y^2,
\\40x^2-2p\pmod{p^2}&\text{if}\  (\f p{13})=1,\ (\f {-2}p)=(\f p{5})=-1\ \&\ p=10x^2+13y^2,
\\0\pmod{p^2}&\text{if}\ (\f{-130}p)=-1,
\end{cases}\end{align*}
where $x$ and $y$ are integers.

{\rm (ii)} We have the congruence
$$\sum_{k=0}^{p-1}(1365k+181)\f{S_k^{(2)}(1/18)}{16^k}\eq p\l(221\l(\f{-2}p\r)-40\l(\f {5}p\r)\r)\pmod{p^2}.$$
Moreover, when $(\f{-10}p)=1$, for any $n\in\Z^+$ the number
$$\f1{(pn)^2}\(\sum_{k=0}^{pn-1}(1365k+181)\f{S_k^{(2)}(1/18)}{16^k}-p\l(\f {-2}p\r)\sum_{k=0}^{n-1}
(1365k+181)\f{S_k^{(2)}(1/18)}{16^k}\)$$
is a $p$-adic integer.
\end{conjecture}
\begin{remark} The imaginary quadratic field $\mathbb Q(\sqrt{-130})$ has class number $4$.
\end{remark}

\section*{Acknowledgment} The author would like to thank the referee for helpful comments.

\end{document}